\documentclass[number,citesort,MSNbibl,seceqn,dvips]{arxbj}
\usepackage{upgreek}

\aid{0}
\volume{18}
\issue{2}
\pubyear{2012}
\firstpage{476}
\lastpage{495}
\doi{10.3150/10-BEJ348}

\makeatletter
\newtheorem{theo}{Theorem}[section]
\newtheorem{lem}[theo]{Lemma}
\newtheorem{prop}[theo]{Proposition}
\newproclaim{defn}[theo]{Definition}
\newremark{ex}[theo]{Example}
\newremark{rem}[theo]{Remark}
\newcommand{\Ph}{\Phi}
\newcommand{\R}{\mathbb{R}}
\newcommand{\rd}{\mathbb {R}^d }
\newcommand{\eqref}[1]{(\ref{#1})}
\makeatother

\begin{document}
\begin{frontmatter}

\title{A class of multivariate infinitely divisible distributions related to arcsine density}
\runtitle{A class of multivariate infinitely divisible
distributions}

\begin{aug}
\author[1]{\fnms{Makoto} \snm{Maejima}\corref{}\thanksref{1}\ead[label=e1]{maejima@math.keio.ac.jp}},
\author[2]{\fnms{V\'{i}ctor} \snm{P\'{e}rez-Abreu}\thanksref{2}\ead[label=e2]{pabreu@cimat.mx}} \and
\author[3]{\fnms{Ken-iti} \snm{Sato}\thanksref{3}\ead[label=e3]{ken-iti.sato@nifty.ne.jp}}
\runauthor{M. Maejima, V. P\'{e}rez-Abreu and K. Sato}
\address[1]{Department of Mathematics, Keio University, 3-14-1,
Hiyoshi, Yokohama 223-8522, Japan.\\\printead{e1}}
\address[2]{Department of Probability and Statistics, Centro de
Investigaci\'{o}n en Matem\'{a}ticas, CIMAT, Apdo. Postal~402, Guanajuato,
Gto. 36000, Mexico. \printead{e2}}
\address[3]{Hachiman-yama 1101-5-103, Tenpaku-ku, Nagoya 468-0074,
Japan.\\ \printead{e3}}
\end{aug}

\received{\smonth{3} \syear{2010}}
\revised{\smonth{11} \syear{2010}}

%
\begin{abstract}
Two transformations $\mathcal{A}_{1}$ and $\mathcal{A}_{2}$ of L\'{e}vy
measures on $\mathbb{R}^{d}$ based on the arcsine density are studied
and their
relation to general Upsilon transformations is considered.
The domains of definition of $\mathcal{A}_{1}$ and $\mathcal{A}_{2}$
are determined and it is
shown that they have the same range.
The class of infinitely divisible distributions on
$\mathbb{R}^{d}$ with L\'{e}vy measures being in the common range
{is} called {the class $A$}
and {any distribution in the class $A$ is} expressed as the law of a
stochastic integral $\int_0^1\cos(2^{-1}\uppi t)\,\mathrm{d}X_t$ with
respect to a L\'{e}vy process $\{X_t\}$.
This new class includes as a proper subclass
the Jurek class of distributions.
It is shown that generalized type $G$
distributions are the image of distributions {in the class $A$}
under a mapping defined
by an appropriate stochastic integral.
$\mathcal{A}_{2}$ is identified as an
Upsilon transformation, while $\mathcal{A}_{1}$ is shown not to be.
\end{abstract}

%

\end{frontmatter}

\section{Introduction}\label{sec1}

Let $I(\mathbb{R }^{d})$ denote the class of all infinitely divisible
distributions on $\mathbb{R }^{d}$.
For $\mu\in I(\mathbb{R}^{d})$, we use the
L\'{e}vy--Khintchine representation of its characteristic function
$\widehat{\mu}(z)$ given by
\begin{eqnarray*}
\label{LevyKhintRep}
\widehat{\mu}(z)&=&\exp\biggl\{ -\frac{1}{2}{}\langle\Sigma z,
z\rangle
 +\mathrm{i} \langle\gamma,z\rangle\\
&&\hphantom{\exp\biggl\{}{} +\int_{\mathbb{R}^{d}} \biggl( \mathrm{e}^{\mathrm{i}\langle
x,z\rangle
} -1- \frac{\mathrm{i}\langle x,z\rangle}{1+ \vert x
\vert^{2}
} \biggr) \nu(\mathrm{d}x) \biggr\}, \qquad   z\in\mathbb{R}^{d},
\end{eqnarray*}
where $\Sigma$ is a symmetric non-negative-definite $d\times d$ matrix,
$\gamma\in\mathbb{R}^{d}$ and $\nu$ is a measure on $\mathbb
{R}^{d}$ (called
the L\'{e}vy measure) satisfying $\nu(\{0\})=0$ and $\int_{\mathbb{R}^{d}
}(1\wedge \vert x \vert^{2})\nu(\mathrm{d}x)<\infty. $
The triplet\vadjust{\goodbreak}
$(\Sigma,\nu,\gamma)$ is called the L\'{e}vy--Khintchine triplet of
$\mu\in
I(\mathbb{R}^{d})$. Let $\mathfrak{M}_{L}(\mathbb{R}^{d})$ denote
the class of
L\'evy measures of $\mu\in I(\mathbb{R}^{d})$.
The class of
$\nu\in\mathfrak{M}_{L}(\mathbb{R}^{d})$ satisfying $\int_{\mathbb
{R}^{d}%
}(1\wedge \vert x \vert)\nu(\mathrm{d}x)<\infty$ is
denoted by
$\mathfrak{M}_{L}^{1}(\mathbb{R}^{d})$.

Let
%
\begin{equation}
a(x;s)=\uppi^{-1}(s-x^{2})^{-1/2} 1_{(-s^{1/2},s^{1/2})}(x), \label{arcs0}
\end{equation}
which is the density of the symmetric arcsine law with parameter $s>0$.
Here $1_{(-s^{1/2},s^{1/2})}(x)$ is the indicator of the interval
$(-s^{1/2},s^{1/2})$.
In \cite{BNPA08}, a symmetric distribution such that its L\'{e}vy measure
has a density $\ell$ of the form 
\[
\ell(x)=\int_{\mathbb{R}_{+}}a(x;s)\rho(\mathrm{d}s), \qquad   x\in
\mathbb{R},
\]
with a measure $\rho$ on $(0,\infty)$ satisfying $\int_{(0,\infty
)}(1\land x)
\rho(\mathrm{d}x)<\infty$ is called a \textit{type $A$ distribution} on
$\mathbb{R}$.
Let $Z$ be a standard normal random variable and $V$ be a positive
infinitely divisible random variable independent of $Z$. The
distribution of
the one-dimensional random variable $V^{1/2}Z$ is infinitely divisible
and is
called type $G$. It is shown in \cite{BNPA08} that an infinitely
divisible distribution $\widetilde{\mu}$ on $\mathbb{R}$ is of type
$G$ if and
only if there exists a type $A$ distribution $\mu$ on $\mathbb{R}$
that gives
a stochastic integral mapping representation
%
\begin{equation}
\widetilde{\mu}=\mathcal{L} \biggl( \int_{0}^{1}( -\log t)%
^{1/2}\,\mathrm{d}X_{t}^{(\mu)} \biggr) . \label{Grep}%
\end{equation}
Here and in what follows, $\mathcal{L}$ means ``the
law of''
and $\{X_{t}^{(\mu)}\}$ means a L\'{e}vy process on~$\mathbb{R}^{d}$ whose
distribution at time 1 is $\mu\in I(\mathbb{R}^{d})$  ($d=1$ in
\eqref{Grep}).

In this paper, we define and study {a class of infinitely divisible}
distributions on $\mathbb{R}^{d}$, {called the class $A$ and
denoted by $A(\rd)$. A distribution in $A(\rd)$ is called a
distribution of
class $A$ in this paper. When $d=1$ and $\mu\in A(\R)$ is symmetric,
$\mu$ is
a type $A$ distribution in \cite{BNPA08}.}
The organization of the paper is the following.

Section~\ref{sec2} introduces the arcsine transformation $\mathcal{A}_{1}$
of L\'{e}vy measures on $\mathbb{R}^{d}$ based on~(\ref{arcs0})
and considers its domain and range. It
is shown that the domain of $\mathcal{A}_{1}$
is~$\mathfrak{M}_{L}^{1}(\mathbb{R}^{d})$. We prove that $\mathcal
{A}_{1}$ is a
one-to-one mapping.
It is shown that the range $\mathfrak{R}(\mathcal{A}_1)$
contains as a proper subclass the class of L\'{e}vy measures of distributions
in the Jurek class~$U(\mathbb{R}^{d})$ studied in \cite{Ju85,MS08}.
The class $U(\mathbb{R}^{d})$ includes several known classes
of multivariate distributions characterized by the radial part of their
L\'{e}vy measures, such as the Goldie--Steutel--Bondesson class
$B(\mathbb{R}%
^{d})$, the class of selfdecomposable distributions $L(\mathbb
{R}^{d})$ and
the Thorin class $T(\mathbb{R}^{d})$; see \cite{BNMS06}. Recently, other
bigger classes than the Jurek class have been discussed in the study of the
extension of selfdecomposability; see \cite{MMS10,Sato10}.

Section~\ref{sec3} deals with the class $A(\mathbb{R}^{d})$ {whose
elements are} defined as infinitely divisible distributions on
$\mathbb{R}^{d}$ {with} L\'{e}vy measures {$\in$}
$\mathfrak{R}(\mathcal{A}_{1})$.
Some probabilistic interpretations of~$A(\mathbb{R}^{d})$
are given and the
relation to the class $G(\mathbb{R}^{d})$ of generalized type $G$
distributions on $\mathbb{R}^{d}$ introduced in \cite{MS08} is
studied. It is
shown that $A(\mathbb{R}^{d})=\Phi_{\cos}(I(\mathbb{R}^{d})),$ where
$\Phi_{\cos}$ is the stochastic integral mapping
%
\begin{equation}\label{Arep}
\Phi_{\cos}(\mu)=\mathcal{L} \biggl( \int_{0}^{1}\cos(2^{-1}\uppi
t)\,\mathrm{d}
X_{t}^{(\mu)} \biggr) , \qquad  \mu\in I(\mathbb{R}^{d}).\vadjust{\goodbreak}
\end{equation}
It is also shown that the class $\mathfrak{M}_L^G (\mathbb{R}^{d})$
is the image of the class $\mathfrak{M}_L^B(\mathbb{R}^{d}%
)\cap\mathfrak{M}_{L}^{1}(\mathbb{R}^{d})$ under $\mathcal{A}_{1}$, where
$\mathfrak{M}_L^G (\mathbb{R}^{d})$ and $\mathfrak{M}_L^B(\mathbb{R}^{d})$
are the classes of L\'{e}vy measures of distributions in $G(\mathbb{R}^{d})$
and $B(\mathbb{R}^{d})$, respectively. In
addition, the class $G(\mathbb{R}^{d})$ is described as the image of
$A(\mathbb{R}^{d})$ under the stochastic integral mapping (\ref{Grep}),
$d\geq1,$ including the multivariate and non-symmetric cases. In order
to prove these facts, a modification $\mathcal{A}_2$ of the transformation
$\mathcal{A}_1$ with the property $\mathfrak{R}(\mathcal{A}_2)=
\mathfrak{R}(\mathcal{A}_1)$ is introduced and utilized effectively.
It is shown that $\mathcal{A}_{2}$ is an Upsilon transformation in
the sense of \cite{BNRT08}. This is in contrast to the fact that
$\mathcal{A}
_{1}$ is not an Upsilon transformation as it is not commuting with a specific
Upsilon transformation, which is different from other cases considered
so far.
Finally, Section~\ref{sec4} contains examples of $\mathcal{A}_{1}$ and
$\mathcal{A}
_{2}$ transformations of L\'{e}vy measures where the modified Bessel function
$K_{0}$ plays an important role.


\section{Arcsine transformation $\mathcal{A}_{1}$ on $\mathbb{R}^{d}$}\label{sec2}


\subsection{Definition and domain}\label{sec2.1}


Besides the arcsine density (\ref{arcs0}), we consider the one-sided arcsine
density
\[
a_{1}(r;s)=2\uppi^{-1}(s-r^{2})^{-1/2} 1_{(0,s^{1/2})}(r) 
\]
with parameter $s>0$.
Then we consider the following arcsine transformation $\mathcal{A}_{1}$
of measures on $\mathbb{R }^{d}$.

\begin{defn}
\label{d1} Let $\nu$ be a measure on $\mathbb{R}^{d}$.
Define the \textit{arcsine transformation} $\mathcal{A}_{1}$ of $\nu$ by
%
\begin{equation}
\mathcal{A}_{1}(\nu)(B)=\int_{\mathbb{R}^{d}\setminus\{0\}}\nu
(\mathrm{d}x)
\int_{0}^{\infty}a_{1}(r;|x|)1_{B} \biggl( r\frac{x}{|x|} \biggr)
\,\mathrm{d}r,  \qquad  B\in\mathcal{B}(\mathbb{R}^{d}), \label{DefArcTrn1}
\end{equation}
for $\nu\in\mathfrak{D}(\mathcal{A}_1)$, where the
domain $\mathfrak{D}(\mathcal{A}_{1})$ is the class of measures $\nu
$ on
$\mathbb{R}^{d}$ such that $\nu(\{0\})=0$ and the right-hand side of
\eqref{DefArcTrn1} belongs to $\mathfrak{M}_{L}(\mathbb{R}^{d})$ as
$B$ runs in $\mathcal{B}(\mathbb{R}^{d})$. The range is
$\mathfrak{R}(\mathcal{A}_{1})=\{\mathcal{A}_{1}(\nu)\dvt\nu\in
\mathfrak{D}(\mathcal{A}_{1})\}$.
\end{defn}

\begin{prop}
\label{DomArcTrn}
$\mathfrak{D}(\mathcal{A}_{1})=\mathfrak{M}_{L}^{1}(\mathbb{R}^{d})$.
\end{prop}

\begin{pf} Suppose that $\nu\in\mathfrak{D}(\mathcal
{A}_{1})$.
Write $\widetilde{\nu}=\mathcal{A}_{1}(\nu)$. Then
%
\begin{equation}
\label{DefArcTrn5}
\int_{\mathbb{R}^{d}}(1\land|x|^{2})\widetilde
{\nu}
(\mathrm{d}x) =c\int_{\mathbb{R}^{d}} \nu(\mathrm{d}x) \int_{0}
^{|x|^{1/2}} (|x| -r^{2})^{-1/2}(1\land r^{2})\,\mathrm{d} r
\end{equation}
with $c=2\uppi^{1/2}$. Changing variables $r=|x|^{1/2} u$ and using
$\int(1\land|x|^{2})\widetilde{\nu}(\mathrm{d}x)<\infty$, we see that
\[
c\int_{\mathbb{R}^{d}} (1\land|x|)\nu(\mathrm{d}x) \int_{0}%
^{1}(1-u^{2})^{-1/2} u^{2}\,\mathrm{d} u<\infty.
\]
Hence $\nu\in\mathfrak{M}_{L}^{1}(\mathbb{R}^{d})$. This shows
$\mathfrak{D}(\mathcal{A}_1)\subset\mathfrak{M}_{L}^{1}(\mathbb{R}^{d})$.\vadjust{\goodbreak}

Next suppose that $\nu\in\mathfrak{M}_{L}%
^{1}(\mathbb{R}^{d})$. Let $\widetilde{\nu}(B)$ denote the
right-hand side
of \eqref{DefArcTrn1}. Then $\widetilde{\nu}$ is a measure on
$\mathbb{R}^{d}$ with $\widetilde{\nu}(\{0\})=0$ and \eqref{DefArcTrn5}
holds. Hence
\[
\int_{\mathbb{R}^{d}}(1\land|x|^{2})\widetilde{\nu}(\mathrm{d}x)
\leq\mathrm{const.} \int_{\mathbb{R}^{d}} (1\land|x|) \nu(\mathrm{d}x)
<\infty.
\]
This shows that $\nu\in\mathfrak{D}(\mathcal{A}_{1})$. Thus
$\mathfrak{M}_{L}^{1}(\mathbb{R}^{d}) \subset\mathfrak{D}(\mathcal{A}_1)$.
\end{pf}


\subsection{One-to-one property}\label{sec2.2}


We next show that the arcsine transformation $\mathcal{A}_{1}$
is one-to-one. Our proof is different from usual proofs of the
one-to-one property by the use of Laplace transform.

Let us prepare a lemma. For a measure $\rho$ on $(0,\infty)$, define
a measure $\sigma_{(\rho)}$ on $(0,\infty)$~by
%
\begin{equation}\label{sigma}
\sigma_{(\rho)}(\mathrm{d}u)= \biggl( \int_{(u,\infty)} \uppi^{-1/2}
(s-u)^{-1/2}\rho(\mathrm{d}s) \biggr) \,\mathrm{d}u.
\end{equation}
This is a fractional integral of order $1/2$.

\begin{lem}
\label{lem2} Let $\rho$ be a measure on $(0,\infty)$ satisfying
%
\begin{equation}
\rho((b,\infty))<\infty  \qquad \mbox{for all $b>0$}. \label{lem2.1}
\end{equation}
Then
%
\begin{equation}
\sigma_{(\sigma_{(\rho)})}(\mathrm{d}u)=\rho((u,\infty))\, \mathrm{d}u,
\label{lem2.2}%
\end{equation}
which implies that $\rho$ is determined by $\sigma_{(\rho)}$ under
condition (\ref{lem2.1}).
\end{lem}

\begin{pf} Using Fubini's theorem, notice that
\begin{eqnarray*}
&& \int_{u}^{\infty}\uppi^{-1/2}(s-u)^{-1/2} \sigma_{(\rho)}(\mathrm
{d}s)\\
&& \quad   =\uppi^{-1}\int_{u}^{\infty}(s-u)^{-1/2}\,\mathrm{d}s\int
_{(s,\infty
)}(v-s)^{-1/2}\rho(\mathrm{d}v)\\
&& \quad   =\uppi^{-1}\int_{(u,\infty)}\rho(\mathrm{d}v)\int_{u}^{v}
(s-u)^{-1/2}(v-s)^{-1/2}\,\mathrm{d}s\\
&& \quad =\rho((u,\infty)),
\end{eqnarray*}
because
\begin{eqnarray*}
\int_{u}^{v}(s-u)^{-1/2}(v-s)^{-1/2}\,\mathrm{d}s&=&\int_{0}^{1}s^{-1/2}
(1-s)^{-1/2}\,\mathrm{d}s\\
&=&B(1/2,1/2)=\uppi,
\end{eqnarray*}
where $B(\cdot,\cdot)$ is the beta function. Hence \eqref{lem2.2} is
true.\vadjust{\goodbreak}~%
\end{pf}
\begin{theo}
\label{t2}
The transformation $\mathcal{A}_{1}$ is one-to-one.
\end{theo}

\begin{pf} Suppose that $\nu,\nu^{\prime}%
\in\mathfrak{M}_{L}^{1}(\mathbb{R}^{d})$ and $\mathcal{A}_{1}(\nu
)=\mathcal{A}_{1} (\nu^{\prime})$. Let $(\lambda, \nu_{\xi})$ and
$(\lambda^{\prime}, \nu^{\prime}_{\xi})$ be polar decompositions
(radial decompositions) of $\nu$ and
$\nu^{\prime}$, respectively, as in Proposition~3.1 of~\cite{Sato10}.
That is, $\lambda$ is a measure on the unit sphere $\mathbb{S}=
\{\xi\in\mathbb{R}^{d}
\dvt  \vert\xi \vert=1\}$ with $0\leq\lambda(\mathbb
{S})\leq
\infty$ and $\nu_{\xi}$, $\xi\in\mathbb{S}$, are
measures on $(0,\infty)$ such that $\nu_{\xi}(E)$ is measurable in
$\xi$ for
each $E\in\mathcal{B}((0,\infty))$, $0<\nu_{\xi}((0,\infty))\leq
\infty$
and
\[
\nu(B)=\int_{\mathbb{S}}\lambda(\mathrm{d}\xi)\int_{(0,\infty)}
1_{B}(u\xi
)\nu_{\xi}(\mathrm{d}u), \qquad   B\in\mathcal{B}(\mathbb{R}^{d});
\]
$\lambda^{\prime}$ and $\nu_{\xi}^{\prime}$ have similar
properties with
respect to $\nu^{\prime}$. It follows from Definition \ref{d1} that
\begin{eqnarray*}
\mathcal{A}_{1}(\nu)(B) & =&\int_{\mathbb{R}^d\setminus\{ 0\}} \nu
(\mathrm{d}x)
\int_0^{|x|^{1/2}} 2\uppi^{-1} (|x|-r^2)^{-1/2} 1_B  \biggl( r\frac{x}{|x|}
 \biggr) \,\mathrm{d}r\\
& =&\int_{\mathbb{S}}\lambda(\mathrm{d}\xi) \int%
_{(0,\infty)}\nu_{\xi}(\mathrm{d}u)\int_0^{u^{1/2}} 2\uppi^{-1}
(u-r^2)^{-1/2}
1_{B}(r\xi)\,\mathrm{d}r\\
& =&\int_{\mathbb{S}}\lambda(\mathrm{d}\xi) \int%
_{0}^{\infty} 1_{B}(r\xi)\,\mathrm{d}r \int_{(r^2,\infty)}2\uppi
^{-1}(u-r^2)^{-1/2}
\nu_{\xi}(\mathrm{d}u)
\end{eqnarray*}
and
\[
\mathcal{A}_{1}(\nu^{\prime})(B)  =\int_{\mathbb{S}} \lambda
^{\prime
}(\mathrm{d}\xi)\int_{0}^{\infty} 1_{B}(r\xi)\,\mathrm{d}r \int
_{(r^2,\infty)}
2\uppi^{-1}(u-r^2)^{-1/2}
\nu_{\xi}^{\prime}(\mathrm{d}u).
\]
These give polar decompositions of $\mathcal{A}_{1}(\nu
)=\mathcal{A}_{1} (\nu^{\prime})$. Hence, by Proposition 3.1 of
\cite{Sato10},
there is a measurable function $c(\xi)$ satisfying $0<c(\xi)<\infty$
such that $\lambda^{\prime}(\mathrm{d}\xi)=c(\xi)\lambda(\mathrm
{d}\xi)$ and,
for $\lambda$-a.e.~$\xi$,
\[
 \biggl( \int_{(r^{2},\infty)} (u-r^{2})^{-1/2}\nu^{\prime}_{\xi}%
(\mathrm{d}u) \biggr) \,\mathrm{d}r=  \biggl( c(\xi)^{-1}\int
_{(r^{2},\infty)}
(u-r^{2})^{-1/2} \nu_{\xi}(\mathrm{d}u) \biggr) \,\mathrm{d}r.
\]
Using a new variable $v=r^{2}$, we see that
\[
 \biggl( \int_{(v,\infty)} (u-v)^{-1/2}\nu^{\prime}_{\xi}(\mathrm
{d}u) \biggr)
\,\mathrm{d}v=  \biggl( c(\xi)^{-1}\int_{(v,\infty)} (u-v)^{-1/2}\nu
_{\xi}
(\mathrm{d}u) \biggr) \,\mathrm{d}v.
\]
Since $\nu_{\xi}$ and $\nu^{\prime}_{\xi}$ satisfy \eqref{lem2.1}
for $\lambda$-a.e. $\xi$, we obtain
$\nu_{\xi}=c(\xi)^{-1} \nu^{\prime}_{\xi}$ for $\lambda$-a.e.
$\xi$ from
Lemma~\ref{lem2}. It follows that $\nu=\nu^{\prime}$.
\end{pf}


\subsection{Range}\label{sec2.3}

Let us show some necessary conditions for $\widetilde{\nu}$ to
belong to the range of $\mathcal{A}_1$.

\begin{prop}
\label{Range2} If $\widetilde{\nu}$ is in $\mathfrak{R}(\mathcal{A}
_{1})$ and not zero measure, then $\widetilde{\nu}$ has a radial decomposition
$(\lambda, \ell_{\xi}(r)\,\mathrm{d}r)$
having the following properties:
\begin{longlist}[(3)]
\item[(1)] $\ell_{\xi}(r)$ is measurable in $(\xi,r)$
and lower semicontinuous in $r\in(0,\infty)$;
\item[(2)] there is $b_{\xi}%
\in(0,\infty]$ such that $\ell_{\xi}(r)>0$ for $r<b_{\xi}$ and,
if $b_{\xi
}<\infty$, then $\ell_{\xi}(r)=0$ for $r\geq b_{\xi}$;
\item[(3)] ${\liminf_{r\downarrow0} \ell_{\xi}(r)>0}$.
\end{longlist}
\end{prop}

\begin{pf} Let $\widetilde{\nu}=\mathcal{A}_{1}(\nu
)$ with
$\nu\in\mathfrak{M}_{L}^{1}(\mathbb{R}^{d})$ and $(\lambda,\nu
_{\xi})$ be a
polar decomposition of $\nu$. Then the proof of Theorem \ref{t2} shows
that $\widetilde{\nu}\in\mathfrak{M}_{L}%
(\mathbb{R}^{d})$ with radial decomposition $(\lambda,\ell(r)\,\mathrm{d}r)$
where
\[
\ell(r)=2\uppi^{-1}\int_{(r^{2},\infty)}(u%
-r^{2})^{-1/2}\nu_{\xi}(\mathrm{d}u).
\]
Then our assertion is proved in the same way as Proposition 2.13 of
\cite{Sato10}.
\end{pf}


\subsection{How big is $\mathfrak{R}(\mathcal{A}_{1})$?}\label{sec2.4}


Several well-known and well-studied classes of multivariate infinitely
divisible distributions are the following. The Jurek class, the class of
selfdecomposable distributions, the Goldie--Steutel--Bondesson class,
the Thorin
class and the class of generalized type $G$ distributions. They are
characterized only by the radial component of their L\'{e}vy measures and
$\Sigma$ and $\gamma$ in the L\'{e}vy--Khintchine triplet are
irrelevant. Among
them, the Jurek class is the biggest. Recently, some classes
bigger than the Jurek class
have been discussed in the study of extension of selfdecomposability (see
\cite{MMS10,Sato10}). Then a natural question is how big
$\mathfrak{R }(\mathcal{A}_{1})$ is. Let $\mathfrak
{M}_{L}^{U}(\mathbb{R}%
^{d})$ be the class of L\'{e}vy measures of distributions in the Jurek class.
The radial component $\nu_{\xi}$ of $\nu\in\mathfrak
{M}_{L}^{U}(\mathbb{R}%
^{d})$ is characterized as
$\nu_{\xi}(\mathrm{d}r)=\ell_{\xi}(r)\,\mathrm{d}r$
with~$\ell_{\xi}(r)$ being measurable in $(\xi,r)$ and decreasing and
right-continuous in $r>0$. We will show below that $\mathfrak{R }%
(\mathcal{A}_{1})$ is strictly bigger than $\mathfrak{M}_{L}%
^{U}(\mathbb{R}^{d})$.

\begin{theo}
\label{Range3}
$\mathfrak{M}_{L}^{U}(\mathbb{R}^{d})\subsetneqq\mathfrak
{R}(\mathcal{A}_{1})$.
\end{theo}

\begin{lem}
\label{lem1a} Let $\rho$ be a $\sigma$-finite measure on $(0,\infty
)$. Then, for
$\alpha>-1$ and $b>0$, the measure $\sigma_{(\rho)}$ in \eqref{sigma}
satisfies
%
\begin{equation}
\label{lem1a-1}\int_{(b,\infty)} u^{\alpha}\sigma_{(\rho
)}(\mathrm{d}u)\leq
C_{1}\int_{(b,\infty)} s^{\alpha+1/2}\rho(\mathrm{d}s)
\end{equation}
and
%
\begin{equation}
\label{lem1a-2}\int_{(0,b]} u^{\alpha}\sigma_{(\rho)}(\mathrm
{d}u)\leq
C_{2} \biggl( \int_{(0,b]} s^{\alpha+1/2}\rho(\mathrm{d}s)+\int
_{(b,\infty)}
s^{-1/2}\rho(\mathrm{d}s) \biggr) ,
\end{equation}
where $C_{1}$ and $C_{2}$ are constants independent of $\rho$.
\end{lem}

\begin{pf}
Let $c=\uppi^{-1/2}$. We have
\[
\int_{(b,\infty)}u^{\alpha}\sigma_{(\rho)}(\mathrm{d}u) =c\int
_{(b,\infty)}
\rho(\mathrm{d}s)\int_{b}^{s}u^{\alpha
}(s-u)^{-1/2}\,\mathrm{d}u
\]
by Fubini's theorem, and
\[
\int_{b}^{s}u^{\alpha}(s-u)^{-1/2}\,\mathrm{d}u=s^{\alpha+1/2}
\int_{b/s}^{1}v^{\alpha}(1-v)^{-1/2}\,\mathrm{d}v \sim
s^{\alpha+1/2}B(\alpha+1,1/2), \qquad   s\rightarrow\infty.
\]
Hence \eqref{lem1a-1} holds. We have
\[
\int_{(0,b]}u^{\alpha}\sigma_{(\rho)}(\mathrm{d}u) =c\int
_{(0,b]}\rho
(\mathrm{d}s)\int_{0}^{s}u^{\alpha}(s-u)^{-1/2}
\,\mathrm{d}u+c\int_{(b,\infty)}\rho(\mathrm{d}s)\int
_{0}^{b}u^{\alpha
}(s-u)^{-1/2}\,\mathrm{d}u.
\]
Notice that
\[
\int_{0}^{s}u^{\alpha}(s-u)^{-1/2}\,\mathrm{d}u=s^{\alpha
+1/2}B(\alpha+1,1/2)
\]
and
\begin{eqnarray*}
  \int_{0}^{b}u^{\alpha}(s-u)^{-1/2}\,\mathrm{d}u&=&s^{-1/2}\int_{0}^{b}
u^{\alpha}(1-u/s)^{-1/2}\,\mathrm{d}u\\
   &\leq& s^{-1/2}\int_{0}^{b}u^{\alpha}(1-u/b)^{-1/2}\,\mathrm{d}
u\\
&=&s^{-1/2}b^{\alpha+1}B(\alpha+1,1/2), \qquad   s>b.
\end{eqnarray*}
Thus \eqref{lem1a-2} holds.
\end{pf}

\begin{pf*}{Proof of Theorem \ref{Range3}}
 Let $\widetilde
{\nu}\in
\mathfrak{M}_{L}^{U}(\mathbb{R}^{d})$ and $c=\int_{\mathbb{R}^d}
(1\wedge|x|^2)\widetilde{\nu}(\mathrm{d}x)$. Then $\widetilde{\nu
}$ has
a polar decomposition $(\lambda,\ell_{\xi}(r)\,\mathrm{d}r)$
mentioned above.
Further, we assume that $\lambda$ is a probability measure and
$\int_{0}^{\infty}(1\wedge r^{2})\ell_{\xi}(r)\,\mathrm{d}r=c$ for
all $\xi$
(see the proof of Proposition 3.1 of \cite{Sato10}, letting
$f(x)=1\wedge
|x|^2$).
Let $\rho_{\xi}$ be a measure on $(0,\infty)$ such that $\rho_{\xi}
((r^{2},\infty))=\ell_{\xi}(r)$ for $r>0$ and let $\eta_{\xi
}=\sigma_{%
(\rho_{\xi})}$. The proof of Lemma \ref{lem2} shows that
\[
\rho_{\xi}((u,\infty))=\int_{(u,\infty)}\uppi
^{-1/2}(s-u)^{-1/2}\eta_{\xi
}(\mathrm{d}s).
\]
Hence we have, for $B\in\mathcal{B}(\mathbb{R}^{d})$,
\[
\widetilde{\nu}(B)=\int_{\mathbb{S}}\lambda(\mathrm{d}\xi)\int
_{0}^{\infty}1_{B}%
(r\xi)\,\mathrm{d}r\int_{(r^{2},\infty)}\uppi
^{-1/2}(s-r^{2})^{-1/2}\eta_{\xi
}(\mathrm{d}s).
\]
We claim that
%
\begin{equation}
\int_{\mathbb{S}}\lambda(\mathrm{d}\xi)\int_{0}^{\infty}(1\wedge
u)\eta_{\xi
}(\mathrm{d}u)<\infty. \label{Range3-2}\
\end{equation}
This will ensure that $\widetilde{\nu}=\mathcal{A}_{1}(\nu)$ for
$\nu\in\mathfrak{D}(\mathcal{A}_1)$ that has polar
decomposition
$(\lambda,\break2^{-1}\uppi^{1/2}\eta_{\xi}(\mathrm{d}r))$. First, notice that
\begin{eqnarray*}
c & =&\frac{1}{2}\int_{0}^{1}u^{1/2}\rho_{\xi}((u,\infty))\,\mathrm
{d}u+\frac
{1}{2}\int_{1}^{\infty}u^{-1/2}\rho_{\xi}((u,\infty))\,\mathrm{d}u\\
& \geq&\frac{1}{3}\rho_{\xi}((1,\infty))+\frac{1}{2}\int
_{1}^{\infty}%
u^{-1/2}\rho_{\xi}((u,\infty))\,\mathrm{d}u.
\end{eqnarray*}
Then, use \eqref{lem1a-1} with $\alpha=0$ to obtain
\begin{eqnarray*}
\int_{(1,\infty)}\eta_{\xi}(\mathrm{d}u)&\leq& C_{1}\int
_{(1,\infty)}
s^{1/2}\rho_{\xi}(\mathrm{d}s)\\
&=&C_{1}\rho_{\xi}((1,\infty))+\frac{C_{1}}{2}\int_{1}^{\infty}%
s^{-1/2}\rho_{\xi}((s,\infty))\,\mathrm{d}s\leq3cC_{1}.
\end{eqnarray*}
Similarly, using \eqref{lem1a-2} with $\alpha=1$,
\begin{eqnarray*}
\int_{(0,1]}u \eta_{\xi}(\mathrm{d}u)&\leq& C_{2} \biggl( \int
_{(0,1]}s^{3/2}
\rho_{\xi}(\mathrm{d}%
s)+\int_{(1,\infty)}s^{-1/2}\rho_{\xi}(\mathrm{d}s) \biggr) \\
&\leq& C_{2} \biggl( \frac{3}{2}\int_{0}^{1}s^{1/2}\rho_{\xi
}((s,1])\,\mathrm{d}s+\int_{(1,\infty)}s^{1/2}\rho_{\xi}(\mathrm
{d}s) \biggr)
\leq5cC_{2}.
\end{eqnarray*}
Hence \eqref{Range3-2} is true. It follows that $\mathfrak{M}_{L}%
^{U}(\mathbb{R}^{d})\subset\mathfrak{R}(\mathcal{A}_{1})$.

To see that the inclusion is strict, consider $\eta\in
\mathfrak{R}(\mathcal{A}_{1})$ defined by
\[
\eta(B) =\int_{\mathbb{S}}\lambda(\mathrm{d}\xi)\int_{0}^{1} 1_{B}(r\xi
)2\uppi^{-1}%
(1-r^{2})^{-1/2}\,\mathrm{d}r.
\]
Then $\eta\notin\mathfrak{M}_{L}^{U}(\mathbb{R}^{d})$, since
$(1-r^{2})^{-1/2}$ is strictly increasing on $(0,1)$.
\end{pf*}


\section{{Distributions of class $A$}}\label{sec3}


\subsection{Definition and representation by another arcsine
transformation $\mathcal{A}_2$}\label{sec3.1}


\begin{defn}
\label{DefTARd} A probability distribution in $I(\mathbb{R}^{d})$ is
said to
be a \textit{distribution {of class~$A$}} on $\mathbb{R }^{d}$ if its
L\'{e}vy measure
${\nu}$ belongs to $\mathfrak{R}(\mathcal{A}_{1})$.
There is no restriction on ${\Sigma}$ and ${\gamma}$ in its
L\'evy--Khintchine triplet. We denote by $A(\mathbb{R}^{d})$ {the
totality of such} distributions on $\mathbb{R}^{d}$.
\end{defn}

Let, for $s>0$,
\[
a_2(r;s)=a_1(r;s^2)=2\uppi^{-1} (s^2-r^2)^{-1/2} 1_{(0,s)}(r).
\]

\begin{defn}
\label{DefA2} Let $\nu$ be a measure on $\mathbb{R}^{d}$.
Define the \textit{arcsine transformation} $\mathcal{A}_{2}$ of $\nu$ by
%
\begin{equation}
\mathcal{A}_{2}(\nu)(B)=\int_{\mathbb{R}^{d}\setminus\{0\}}\nu
(\mathrm{d}x)
\int_{0}^{\infty}a_{2}(r;|x|)1_{B} \biggl( r\frac{x}{|x|} \biggr)
\,\mathrm{d}r,  \qquad  B\in\mathcal{B}(\mathbb{R}^{d}). \label
{DefArcTrn2}%
\end{equation}
The domain $\mathfrak{D}(\mathcal{A}_{2})$ is defined to be the class of
$\nu$ such that $\nu(\{0\})=0$ and the right-hand side of
\eqref{DefArcTrn2} belongs to
$\mathfrak{M}_{L}(\mathbb{R}^{d})$.\vspace*{-1pt}
\end{defn}

For any measure $\nu$ on $\mathbb{R}^{d}$ with $\nu(\{0\})=0$,
define a
measure $\nu^{(2)}$ on $\mathbb{R}^{d}$ by
\[
\nu^{(2)}(B)=\int_{\mathbb{R}^{d}\setminus\{0\}} 1_B \biggl( |x|^2
\frac{x}{|x|} \biggr)\nu(\mathrm{d}x), \qquad   B\in\mathcal
{B}(\mathbb{R}^d).
\]
The transformation $\nu\mapsto\nu^{(2)}$ is one-to-one, since we have
$\nu=(\nu^{(2)})^{(1/2)}$, defining $\rho^{(1/2)}$ for $\rho$ with
$\rho(\{0\})=0$ as
\[
\rho^{(1/2)}(B)=\int_{\mathbb{R}^{d}\setminus\{0\}} 1_B \biggl( |x|^{1/2}
\frac{x}{|x|} \biggr)\rho(\mathrm{d}x).
\]
The following propositions give the connections between $\mathcal{A}_1$
and $\mathcal{A}_2$.\vspace*{-1pt}

\begin{prop}\label{conn1}
$\mathfrak{D}(\mathcal{A}_2)=\mathfrak{M}_L(\mathbb{R}^d)$.\vspace*{-1pt}
\end{prop}

\begin{prop}\label{conn2}
$\nu\in\mathfrak{M}_{L}(\mathbb{R}^{d})$ if and only if
$\nu^{(2)}\in\mathfrak{M}_{L}^{1}(\mathbb{R}^{d})$, and in this case
\[
\label{p-trRep-1}
\mathcal{A}_{2}(\nu)=\mathcal{A}_{1}\bigl(\nu^{(2)}\bigr).\vspace*{-1pt}
\]
\end{prop}

\begin{prop}\label{conn3}
$\mathfrak{R}(\mathcal{A}_1)=\mathfrak{R}(\mathcal{A}_2)$.\vspace*{-1pt}
\end{prop}

\begin{pf*}{Proof of Propositions \ref{conn1}--\ref{conn3}} Since
$\int_{\mathbb{R}^d}f(x)\nu^{(2)}(\mathrm{d}x)=\int_{\mathbb{R}^d}
f(|x|x)\nu(\mathrm{d}x)$ for any non-negative measurable function $f$,
we have $\int_{\mathbb{R}^d}(1\wedge|x|)\nu^{(2)}(\mathrm{d}x)=
\int_{\mathbb{R}^d}(1\wedge|x|^2)\nu(\mathrm{d}x)$ and the first two
propositions follow. Proposition \ref{conn3} is a
direct consequence of them.\vspace*{-1pt}
\end{pf*}

We will also use the following fact.\vspace*{-1pt}

\begin{prop}\label{conn4}
$\mathcal{A}_2$ is a one-to-one transformation.\vspace*{-1pt}
\end{prop}

\begin{pf} This follows from Theorem \ref{t2} and
Proposition \ref{conn2}.\vspace*{-1pt}
\end{pf}

The usefulness of the introduction of $\mathcal{A}_2$ is based on the
next property.\vspace*{-1pt}

\subsection{Transformation $\mathcal{A}_{2}$ as an Upsilon transformation}\label{sec3.2}\vspace*{-1pt}


Barndorff-Nielsen, Rosi\'{n}ski and
Thorbj\o rnsen \cite{BNRT08} considered
general Upsilon transformations (see \cite{BNM08,Sa07}).
Given a measure $\tau$ on $(0,\infty)$, a transformation $\Upsilon
_{\tau}$
from measures on\vadjust{\goodbreak} $\mathbb{R}^{d}$ into $\mathfrak{M}_{L}(\mathbb
{R}^{d})$ is
called an Upsilon transformation associated to $\tau$ (or with dilation
measure $\tau$) when
%
\begin{equation}
\Upsilon_{\tau}(\nu)(B)=\int_{0}^{\infty}\nu(u^{-1}B)\tau
(\mathrm{d}u), \qquad
B\in\mathcal{B}(\mathbb{R}^{d}). \label{GeUpsMap}%
\end{equation}
The domain of $\Upsilon_{\tau}$ is the class of $\sigma$-finite
measures $\nu$
such that the right-hand side of~\eqref{GeUpsMap} is a measure in
$\mathfrak{M}_{L}(\mathbb{R}^{d})$.

\begin{theo}
\label{UpsIdent} $\mathcal{A}_2$ is an Upsilon transformation given by
%
\begin{equation}
\mathcal{A}_{2}(\nu)(B)=\int_{0}^{1} \nu(u^{-1}B) 2\uppi^{-1}
(1-u^{2})^{-1/2}\,\mathrm{d}u,  \qquad   B\in\mathcal{B}(\mathbb{R}^{d})
\label{PrUpsIdent1}%
\end{equation}
for $\nu\in\mathfrak{M}_{L}(\mathbb{R}^{d})$.
\end{theo}

\begin{pf} Let $(\lambda,\nu_{\xi})$ be a polar
decomposition of
$\nu\in\mathfrak{M}_{L}(\mathbb{R}^{d})$. Then, with $c=2\uppi^{-1}$,
\begin{eqnarray*}
\mathcal{A}_{2}(\nu)(B) & =&c\int_{\mathbb{S}}\lambda(\mathrm
{d}\xi) \int%
_{0}^{\infty}1_{B}(r\xi)\,\mathrm{d}r\int_{(r,\infty)}(s^{2}-r^{2}%
)^{-1/2}\nu_{\xi}(\mathrm{d}s)\\
& =&c\int_{\mathbb{S}}\lambda(\mathrm{d}\xi)\int_{0}^{\infty}\nu
_{\xi
}(\mathrm{d}s) \int_{0}^{s}1_{B}(r\xi)(s^{2}-r^{2})^{-1/2}\,\mathrm
{d}r\\
& =&c\int_{\mathbb{S}}\lambda(\mathrm{d}\xi)\int_{0}^{\infty}\nu
_{\xi
}(\mathrm{d}s) \int_{0}^{1}1_{B}(us\xi)(1-u^{2})^{-1/2}\,\mathrm{d}u\\
& =&c\int_{0}^{1}(1-u^{2})^{-1/2}\,\mathrm{d}u\int_{\mathbb{S}}\lambda
(\mathrm{d}\xi)\int_{0}^{\infty}1_{B}(us\xi)\nu_{\xi}(\mathrm
{d}s)\\
& =&c\int_{0}^{1}(1-u^{2})^{-1/2}\,\mathrm{d}u\int_{\mathbb
{R}^{d}}1_{B}%
(ux)\nu(\mathrm{d}x),
\end{eqnarray*}
which shows \eqref{PrUpsIdent1}.

Conversely, let $\nu$ be a $\sigma$-finite measure on $\mathbb{R}^{d}$,
define $\widetilde{\nu}(B)$ by the right-hand side of~\eqref{PrUpsIdent1},
and suppose that $\widetilde{\nu}\in\mathfrak{M}_{L}(\mathbb{R}^{d})$.
Then $0=\widetilde{\nu}(\{0\})=\nu(\{0\})$ and the calculation above shows
that $\nu\in\mathfrak{D}(\mathcal{A}_2)$ and $\mathcal{A}_2(\nu)=
\widetilde{\nu}$. Hence $\mathcal{A}_2=\Upsilon_{\tau}$ of the form
\eqref{PrUpsIdent1}.
\end{pf}

The mapping $\mathcal{A}_{1}$ is not an Upsilon transformation for any
dilation measure $\tau$. This remarkable result will be proved in
Section~\ref{sec3.8},
as a byproduct of Theorem \ref{noncom} shown in Section~\ref{sec3.7}.

\subsection{Stochastic integral representation of $A(\mathbb{R}^{d})$}\label{sec3.3}


We study a probabilistic interpretation of
distributions of {class $A$}, representing them by stochastic integrals
with respect to
L\'{e}vy processes.

Let $T\in(0,\infty)$ and let $f(t)$ be a square-integrable function on
$[0,T]$. Then the stochastic integral $\int_{0}^{T}f(t)\,\mathrm
{d}X_{t}^{(\mu
)}$ is\vadjust{\goodbreak} defined for any $\mu\in I(\mathbb{R}^{d})$ and is infinitely divisible.
Define the stochastic integral mapping $\Phi_{f}$ based on $f$ as
$\Phi_{f}(\mu)=\mathcal{L} ( \int_{0}^{T}f(t)\,\mathrm
{d}X_{t}^{(\mu
)} )$ for $\mu\in I(\mathbb{R}^{d})$.
If $\mu\in I(\mathbb{R}^{d})$ has the L\'evy--Khintchine triplet
$(\Sigma
,\nu,\gamma)$, then, as in \cite{Sa06a,Sa06b,Sato10}, $\widetilde
{\mu}
=\Phi_{f}(\mu)$ has the
triplet $(\widetilde{\Sigma},\widetilde{\nu},\widetilde{\gamma})$
expressed
as
%
\begin{eqnarray}\label{stochint2}
\widetilde{\Sigma}&=&\int_{0}^{T}f(t)^{2} \Sigma \,\mathrm
{d}t,%
\\[-2pt]\label{stochint3}
\widetilde{\nu}(B)&=&\int_{0}^{T}\mathrm{d}t\int_{\mathbb{R}^{d}}1_{B}
(f(t)x) \nu(\mathrm{d}x), \qquad   B\in\mathcal{B}(\mathbb{R}^{d}%
),\\[-2pt]\label{stochint4}
\widetilde{\gamma}&=&\int_{0}^{T}f(t)\,\mathrm{d}t \biggl( \gamma+\int%
_{\mathbb{R}^{d}}x  \biggl( \frac{1}{1+|f(t)x|^{2}}-\frac
{1}{1+|x|^{2}} \biggr)
\nu(\mathrm{d}x) \biggr) . %
\end{eqnarray}

Let us characterize the class $A(\mathbb{R}^{d})$ as the range of a stochastic
integral mapping.

\begin{theo}
\label{StIntRep} Define $\Ph_{\cos}$ by \eqref{Arep}.
Then $\Phi_{\cos}$ is a one-to-one mapping and
%
\begin{equation}
A(\mathbb{R}^{d})=\Phi_{\cos}(I(\mathbb{R}^{d})). \label{StIntRep2}
\end{equation}
\end{theo}

\begin{pf} Let $\widetilde{\mu}\in A(\mathbb{R}^{d})$
with triplet $(\widetilde{\Sigma},\widetilde{\nu
},\widetilde{\gamma})$. Then $\widetilde{\nu}\in\mathfrak
{R}(\mathcal{A}_{2})$
by Proposition \ref{conn3} and hence we have \eqref{PrUpsIdent1}
with some $\nu\in\mathfrak{M}_{L}(\mathbb{R}^{d})$ by Theorem \ref
{UpsIdent}.
Let $t=g(u)=\int%
_{u}^{1}2\uppi^{-1}(1-v^{2})^{-1/2}\,\mathrm{d}v= 2\uppi^{-1}\arccos
(u)$ for
$0<u<1$. Then $u=\cos(2^{-1}\uppi t)$ for $0<t<1$. Thus
\[
\widetilde{\nu}(B)=-\int_{0}^{1}\mathrm{d}g(u)\int_{\mathbb{R}^{d}}
1_{B}(ux)\nu(\mathrm{d}x)=\int_{0}^{1}\mathrm{d}t\int_{\mathbb{R}^{d}}
1_{B}(x\cos(2^{-1}\uppi t))\nu(\mathrm{d}x).
\]
That is, \eqref{stochint3} is satisfied with $T=1$ and $f(t)=\cos
(2^{-1}\uppi
t)$. Using $\nu$, $\widetilde{\Sigma}$ and $\widetilde{\gamma}$,
we can find $\Sigma$ and $\gamma$ satisfying
\eqref{stochint2} and \eqref{stochint4}. Let $\mu$ be the
distribution in
$I(\mathbb{R}^{d})$ with triplet $(\Sigma,\nu,\gamma)$.
Then $\widetilde{\mu}=\Phi_{\cos}(\mu)$. Hence $A(\mathbb
{R}^{d})\subset
\Phi_{\cos}(I(\mathbb{R}^{d}))$.

Conversely, suppose that $\widetilde{\mu}\in\Phi_{\cos}(I(\mathbb
{R}^{d}))$.
Then $\widetilde{\mu}=\Phi_{\cos}(\mu)$ for some $\mu\in
I(\mathbb{R}^{d})$
by a~similar argument, we can show that $\Phi_{\cos}%
(I(\mathbb{R}^{d}))\subset A(\mathbb{R}^{d})$.

The mapping $\Phi_{\cos}$ is one-to-one, since $\nu$ is determined by
$\widetilde{\nu}$ (Proposition \ref{conn4}) and $\Sigma$ and
$\gamma$ are
determined by $\widetilde{\Sigma}$, $\widetilde{\gamma}$ and $\nu
$.\vspace*{-3pt}
\end{pf}

\subsection{\texorpdfstring{$\Upsilon^{0}$ transformation}{Upsilon^{0} transformation}}\label{sec3.4}\vspace*{-3pt}


We use $\Upsilon$ and $\Upsilon^{0}$ defined by
\begin{eqnarray*}
\Upsilon(\mu)&=&\mathcal{L} \biggl( \int_0^1(-\log t)\,\mathrm
{d}X_t^{(\mu)}
 \biggr), \qquad  \mu\in I(\mathbb{R}^{d}),\\[-2pt]
\Upsilon^{0}(\nu)(B)&=&\int_{0}^{\infty}\nu(u^{-1}%
B)\mathrm{e}^{-u}\,\mathrm{d}u,  \qquad   B\in\mathcal{B}(\mathbb{R}^{d}).
\end{eqnarray*}
Let $\mathfrak{M}_{L}^{B}(\mathbb{R }^{d})$ be the class of L\'evy
measures of
the Goldie--Steutel--Bondesson class $B(\mathbb{R }^{d})$. In \cite
{BNMS06}, it
is shown that $\Upsilon(I(\mathbb{R}%
^{d}))=B(\mathbb{R}^{d})$ and $\Upsilon^{0}(\mathfrak{M}_{L}(\mathbb
{R }^{d}))
=\mathfrak{M}%
_{L}^{B}(\mathbb{R }^{d})$. Both $\Upsilon^{0}$ and
$\Upsilon$ are one-to-one.\vadjust{\goodbreak}

\begin{prop}
\label{Ups0} Let $\nu\in\mathfrak{M}_{L}(\mathbb{R}^{d})$. Then
$\Upsilon
^{0}(\nu)\in\mathfrak{M}_{L}^{1}(\mathbb{R}^{d})$ if and only if
$\nu
\in\mathfrak{M}_{L}^{1}(\mathbb{R}^{d})$.
\end{prop}

\begin{pf} Notice that
\begin{eqnarray*}
\int_{|x|\leq1} |x|\Upsilon^{0}(\nu)(\mathrm{d}x)&=&
\int_{\mathbb{R}^{d}}|x|\nu(\mathrm{d}x) \int_{0}%
^{1/|x|}u\mathrm{e}^{-u} \,\mathrm{d}u,\\
\int_{|x|>1}\Upsilon^{0}(\nu)(\mathrm{d}x)&=&\int_%
{\mathbb{R}^{d}}\mathrm{e}^{-1/|x|}\nu(\mathrm{d}x),
\end{eqnarray*}
to see the equivalence.\vspace*{2pt}
\end{pf}

\subsection{A representation of completely monotone functions}\label{sec3.5}\vspace*{2pt}


In \cite{MS08}, the class of generalized type $G$ distributions on
$\mathbb{R}^{d}$, denoted by $G(\mathbb{R}^{d})$, is defined as follows:
$\mu\in G(\mathbb{R}^{d})$ if and only if the radial component $\nu
_{\xi}$ of
the L\'{e}vy measure of $\mu$ satisfies $\nu_{\xi}(\mathrm
{d}r)=g_{\xi}%
(r^{2})\,\mathrm{d}r$, where $g_{\xi}(u)$ is measurable in $(\xi,u)$ and
completely monotone in $u>0$.
$\mathfrak{M}_{L}^{G}(\mathbb{R}^{d})$ denotes the class of all
L\'{e}vy measures of $\mu\in G(\mathbb{R}^{d})$. We use the following result
when dealing with $G(\mathbb{R}^{d})$. It is a result on the arcsine
transformation representation of a function $g(r^{2})$ when $g$ is completely
monotone on $(0,\infty)$.\looseness=1

\begin{prop}
\label{Repgr2} Let $g(u)$ be a function on $(0,\infty
)$. Then the following three conditions are equivalent.
\begin{longlist}[(b)]
\item[(a)]
The function $g(u)$ is completely monotone on $(0,\infty)$ and satisfies
%
\begin{equation}
\label{Repgr2-1}\int_{0}^{\infty}(1\land r^{2}) g(r^{2})\,\mathrm
{d}r<\infty.
\end{equation}
\item[(b)] There exists a completely monotone function $h(s)$ on
$(0,\infty)$ satisfying
%
\begin{equation}
\label{Repgr2-2}\int_{0}^{\infty}(1\land s) h(s)\,\mathrm{d}s<\infty
\end{equation}
such that
%
\begin{equation}\label{Repgr2-3}
g(r^{2})=\int_{0}^{\infty} a_{1}(r;s)h(s)\,\mathrm{d}s,  \qquad  r>0.
\end{equation}
\item[(c)] There exists a measure $\rho$ on $(0,\infty)$ satisfying
\[
\label{Repgr2-4}\int_{0}^{\infty}(1\land s) \rho(\mathrm
{d}s)<\infty
\]
such that
%
\begin{equation}
\label{Repgr2-5}g(r^{2})=\int_{(0,\infty)} a_1(r;s)\Upsilon
^{0}(\rho)
(\mathrm{d}s), \qquad   r>0.
\end{equation}
\end{longlist}
\end{prop}

\begin{pf} (a)${}\Rightarrow{}$(b): From Bernstein's theorem,
there exists a measure $Q$ on $[0,\infty)$ such that
%
\begin{equation}
g(u)=\int_{[0,\infty)}\mathrm{e}^{-uv}Q(\mathrm{d}v), \qquad   u>0.
\label{Repgr2-6}%
\end{equation}
It follows from \eqref{Repgr2-1} that $Q(\{0\})=0$, since $Q(\{0\}
)=\lim
_{u\rightarrow\infty}g(u)$. We need the fact that the one-dimensional Gaussian
density $\varphi(x;t)$ of mean 0 and variance $t$ is the arcsine
transform of
the exponential distribution with mean $t>0$. More precisely,
%
\begin{equation}
\varphi(x;t)={(2\uppi t)}^{-1/2}\mathrm{e}^{-x^{2}/(2t)}={t}^{-1} \int%
_{0}^{\infty}\mathrm{e}^{-s/t}a(x;2s)\,\mathrm{d}s, \qquad  t>0, x\in
\mathbb{R}.
\label{RGausAT}%
\end{equation}
This is the well-known Box--Muller method to generate normal random variables.
Its proof can be given by change of variables $s=tu+x^2/2$.
Using \eqref{RGausAT}, we have
\begin{eqnarray*}
g(r^{2}) & =&\int_{(0,\infty)}\mathrm{e}^{-r^{2}v}Q(\mathrm{d}v)
=\int_{(0,\infty)}v^{1/2}Q(\mathrm{d}v)\int_{r^{2}/2}^{\infty}
\mathrm{e}^{-2sv}2\uppi^{-1/2}(2s-r^{2})^{-1/2}\,\mathrm{d}s\\
& =&\int_{r^{2}}^{\infty}\uppi^{-1/2}(s-r^{2})^{-1/2}\,\mathrm{d}s \int%
_{(0,\infty)}\mathrm{e}^{-sv}v^{1/2}Q(\mathrm{d}v)
=\int_{0}^{\infty}a_{1}(r;s)h(s)\,\mathrm{d}s,
\end{eqnarray*}
where
%
\begin{equation}
h(s)=2^{-1}\uppi^{1/2}\int_{(0,\infty)}\mathrm
{e}^{-sv}v^{1/2}Q(\mathrm{d}v).
\label{Repgr2-7}%
\end{equation}
Applying Proposition \ref{DomArcTrn} for $d=1$, we see \eqref{Repgr2-2}
from \eqref{Repgr2-1}.

(b)${}\Rightarrow{}$(c): There is $\rho\in\mathfrak{M}_{L}(\mathbb
{R})$ such that
$h(s)\,\mathrm{d}s= \Upsilon^{0}(\rho)$ (see Theorem A of \cite
{BNMS06}). Since~$\Upsilon^{0}(\rho)$ is concentrated on $(0,\infty)$, $\rho$ is concentrated
on $(0,\infty)$. Using Proposition~\ref{Ups0}, we see that $\int%
_{(0,1]}s \rho(\mathrm{d}s)<\infty$.

(c)${}\Rightarrow{}$(a): It follows from Proposition \ref{Ups0} that
$\int_{(0,\infty)}(1\wedge s) \Upsilon^{0}(\rho)(\mathrm
{d}s)<\infty$. Hence it follows
from \eqref{Repgr2-5} that $g(r^{2})$ satisfies \eqref{Repgr2-1} (use
Proposition
\ref{DomArcTrn} for $d=1$). Finally let us prove that $g(u)$ is completely
monotone. There is a completely monotone function~$h(s)$ such that
$\Upsilon^{0}(\rho)(\mathrm{d}s)=h(s)\,\mathrm{d}s$ (see Theorem A of
\cite{BNMS06} again). Hence we can find a~measure~$R$ on $[0,\infty)$
such that
$h(s)=\int_{[0,\infty)}\mathrm{e}^{-sv}R(\mathrm{d}v)$, $s>0$.
We have $R(\{0\})=0$. Thus
\[
g(r^{2})=\int_{r^{2}}%
^{\infty}2\uppi^{-1}(s-r^{2})^{-1/2}\,\mathrm{d}s\int_{(0,\infty)}
\mathrm{e}%
^{-sv}R(\mathrm{d}v)=\int_{(0,\infty)}\mathrm{e}^{-r^{2}v}2\uppi^{-1/2}
v^{-1/2}R(\mathrm{d}v),
\]
where the last equality is from the same calculus as in the proof that
(a)${}\Rightarrow{}$(b). Now we see that $g(u)$ is completely
monotone.
\end{pf}

\subsection{A representation of $G(\mathbb{R}^{d})$ in terms of
$\mathcal{A}%
_{1}$}\label{sec3.6}


We now give an alternative representation for L\'{e}vy measures of
distributions in $G(\mathbb{R }^{d})$.

\begin{theo}
\label{GenTGRep} Let $\widetilde{\mu}$ be an infinitely divisible
distribution
on $\mathbb{R}^{d}$ with the L\'evy--Khintchine triplet $(\widetilde
{\Sigma
},\widetilde{\nu},\widetilde{\gamma})$. Then the following three conditions
are equivalent.
\begin{longlist}[(b)]
\item[(a)]     $\widetilde{\mu}\in
G(\mathbb{R}^{d})$.
\item[(b)]  $\widetilde{\nu
}=\mathcal{A}%
_{1}(\nu)$ with some $\nu\in\mathfrak{M}_{L}^{B}(\mathbb
{R}^{d})\cap
\mathfrak{M}_{L}^{1} (\mathbb{R}^{d}) $.
\item[(c)]
$\widetilde{\nu}=\mathcal{A}_{1}(\Upsilon^{0}(\rho))$ with some
$\rho
\in\mathfrak{M}_{L}^{1}(\mathbb{R}^{d})$.
\end{longlist}
  In condition
\textup{(b)}
or  \textup{(c)}, the representation of $\widetilde{\nu}$ by $\nu$
or $\rho$
is unique.
\end{theo}

\begin{pf} (a)${}\Rightarrow{}$(b): By definition of
$G(\mathbb{R}^{d})$ and Proposition \ref{Repgr2}, the L\'evy measure
$\widetilde{\nu}$ of~$\widetilde{\mu}$
has polar decomposition $(\lambda, g_{\xi}(r^{2})\,\mathrm{d}r)$ where
$g_{\xi
}(u)$ is measurable in $(\xi,u)$ and\vspace*{2pt} satisfies~\eqref{Repgr2-2} and
\eqref{Repgr2-3} with $g_{\xi}(r^{2})$ and $h_{\xi}(s)$ in place of
$g(r^2)$ and $h(s)$.
The measure~$Q_{\xi}$ in the representation \eqref{Repgr2-6} of
$g_{\xi}(u)$
has the property that $Q_{\xi}(E)$ is measurable in $\xi$ for every
Borel set
$E$ in $[0,\infty)$ (see Remark 3.2 of \cite{BNMS06}). Hence the function
$h_{\xi}(s)$ given as in \eqref{Repgr2-7} is measurable in $(\xi
,s)$. Thus
we have $\widetilde{\nu}=\mathcal{A}_{1}(\nu)$, letting $\nu$
denote the
L\'evy measure with polar decomposition $(\lambda, h_{\xi
}(s)\,\mathrm{d}s)$.
Notice that $\nu\in\mathfrak{M}_{L}^{B} (\mathbb{R}^{d})$ and show that
$\nu\in\mathfrak{M}_{L}^{1}(\mathbb{R}^{d})$ by an argument similar to
the proof of Proposition \ref{DomArcTrn}.

(b)${}\Rightarrow{}$(c): Use Proposition \ref{Ups0} and the representation
of $\mathfrak{M}_{L}^{B}(\mathbb{R}^{d})$ by $\Upsilon^0$.

(c)${}\Rightarrow{}$(a): Use Proposition \ref{Ups0} again together with
Proposition \ref{Repgr2}.

Uniqueness of the representations comes from Theorem \ref{t2} and that
of $\Upsilon^0$.\vspace*{-3pt}
\end{pf}

\subsection{$G(\mathbb{R}^{d})$ as image of $A(\mathbb{R}^{d})$
under a
stochastic integral mapping}\label{sec3.7}\vspace*{-3pt}


Following \cite{MN08}, we define the transformation $\Upsilon_{\alpha
,\beta
}(\nu)$ for $-\infty<\alpha<2$ and $0<\beta\leq2$. For a measure
$\nu$ on
$\mathbb{R}^{d}$ with $\nu(\{0\})=0$ define\vspace*{-1pt}
\[
\Upsilon_{\alpha,\beta}(\nu)(B)=\int_{0}^{\infty}\nu
(s^{-1}B)\beta
s^{-\alpha-1}\mathrm{e}^{-s^{\beta}}\,\mathrm{d}s, \qquad   B\in\mathcal{B}
(\mathbb{R}^{d}),\label{UpsMapalbe}%
\]
whenever the right-hand side gives a measure in $\mathfrak{M}_{L}%
(\mathbb{R}^{d})$. This definition is different from that of \cite
{MN08} in
the constant factor $\beta$. A special case with $\beta=1$ coincides
with the
transformation of L\'{e}vy measures in the stochastic integral mapping
$\Psi_{\alpha}$ studied by Sato~\cite{Sa06b}. Of particular interest
in this
work is the mapping $\Upsilon_{-2,2}$. Notice that $\Upsilon_{-1,1}%
=\Upsilon^{0}$.\vspace*{-1pt}

\begin{prop}
\label{dom-2,2} $\mathfrak{D}(\Upsilon_{-2,2})=\mathfrak{M}_{L}
(\mathbb{R}^{d})$. The mapping $\Upsilon_{-2,2}$ is one-to-one.\vspace*{-1pt}
\end{prop}

\begin{pf} Let $\widetilde{\nu}(B)=\int_{0}^{\infty
}\nu
(s^{-1}B)2s\mathrm{e}^{-s^{2}}\,\mathrm{d}s$. Then\vspace*{-1pt}
\begin{eqnarray*}
\int_{\mathbb{R}^{d}}(1\wedge|x|^{2})\widetilde{\nu}(\mathrm
{d}x)&=&\int%
_{0}^{\infty}2s\mathrm{e}^{-s^{2}}\,\mathrm{d}s\int_{\mathbb{R}^{d}}
(1\wedge|sx|^{2})\nu(\mathrm{d}x)\\[-2pt]
&=&\int_{\mathbb{R}^{d}}|x|^{2}\nu(\mathrm{d}x)\int_{0}^{1/|x|}
2s^{3}\mathrm{e}^{-s^{2}}\,\mathrm{d}s+\int_{\mathbb{R}^{d}}\nu
(\mathrm{d}%
x)\int_{1/|x|}^{\infty} 2s\mathrm{e}^{-s^{2}}\,\mathrm{d}s.
\end{eqnarray*}
Observe that $\int_{0}^{1/|x|}2s^{3}\mathrm{e}^{-s^{2}}\,\mathrm{d}s$ is
convergent as $|x|\downarrow0$ and $\sim2^{-1}|x|^{-4}$ as
$|x|\rightarrow
\infty$ and $\int_{1/|x|}^{\infty}2s\mathrm{e}^{-s^{2}}\,\mathrm
{d}s$ is
$\sim\mathrm{e}^{-1/|x|^{2}}$ as $|x|\downarrow0$ and convergent\vadjust{\goodbreak} as
$|x|\rightarrow\infty$. Then we see that $\int_{\mathbb
{R}^{d}}(1\wedge
|x|^{2})\widetilde{\nu}(\mathrm{d}x)$ is finite if and only if $\int
_{\mathbb{R}^{d}}(1\wedge|x|^{2})\nu(\mathrm{d}x)$ is finite. To
prove that
$\Upsilon_{-2,2}$ is one-to-one, make an argument similar to the proof of
Proposition 4.1 of \cite{Sa06b}.\vspace*{-2pt}
\end{pf}

  The following result is needed in giving the
characterization of
$G(\mathbb{R }^{d})$ in terms of distributions of class $A$. It
shows that $\mathcal{A}_{1}$ and $\Upsilon^{0}$ are not commutative. However,
$\mathcal{A}_{2}$ and~$\Upsilon^{0}$ are commutative, both being Upsilon
transformations with  domain  $\mathfrak{M}_{L}(\mathbb{R}^{d})$.\vspace*{-2pt}

\begin{theo}
\label{noncom} It holds that
\[
\label{noncom-1}
\Upsilon_{-2,2}(\mathcal{A}_{1}(\rho)) =\mathcal{A}
_{1}(\Upsilon^{0}(\rho))  \qquad \mbox{for }\rho\in\mathfrak{M}_{L}^{1}
(\mathbb{R}^{d}).\vspace*{-2pt}
\]
\end{theo}

\begin{pf} Let $\rho\in\mathfrak{M}_{L}^{1}%
(\mathbb{R}^{d})$, $\nu=\mathcal{A}_{1}(\rho)$ and $\widetilde{\nu}=
\Upsilon_{-2,2}(\nu)$ with polar decompositions $(\lambda,\rho_{\xi})$,
$(\lambda,\nu_{\xi})$ and $(\lambda,\widetilde{\nu}_{\xi})$,
respectively.
From Lemma 2.5 of \cite{MN08} we have $\widetilde{\nu}_{\xi
}(\mathrm{d}r)
=rg_{\xi}(r^{2})\,\mathrm{d}r$ with $rg_{\xi}(r^{2})=2r\int
_{0}^{\infty}s^{-2}
\mathrm{e}^{-r^{2}/s^{2}}\nu_{\xi}(\mathrm{d}s)$. Hence
\begin{eqnarray*}
r{g}_{\xi}(r^{2}) & =&2r\int_{0}^{\infty}\mathrm{e}^{-r^{2}/s^{2}}
s^{-2}\,\mathrm{d}s\int_{(0,\infty)}a_{1}(s;u)\rho_{\xi}(\mathrm
{d}u)\\[-2pt]
& =&\int_{0}^{\infty}\mathrm{e}^{-t}t^{-1/2}\,\mathrm{d}t\int
_{(0,\infty)}
a_{1}(t^{-1/2}r;u)\rho_{\xi}(\mathrm{d}u)\\[-2pt]
& =&\int_{0}^{\infty}\mathrm{e}^{-t}\,\mathrm{d}t\int_{0}^{\infty}
a_{1}(r;tu)\rho_{\xi}(\mathrm{d}u)\\[-2pt]
& =&\int_{0}^{\infty}a_{1}(r;u)\Upsilon^{0}(\rho_{\xi})(\mathrm{d}u).
\end{eqnarray*}
It follows that $\widetilde{\nu}=\mathcal{A}_{1}
(\Upsilon^{0}(\rho))$.\vspace*{-2pt}
\end{pf}

The following result shows that $G(\mathbb{R}^{d})$ is the class of
distributions of stochastic integrals with respect L\'{e}vy processes with
distribution {of class $A$} at time $1$. This is a multivariate and not
necessarily
symmetric generalization of \eqref{Grep}.\vspace*{-2pt}

\begin{theo}
\label{noncom1} Let
\[
\Psi_{-2,2}(\mu)=\mathcal{L} \biggl( \int_{0}^{1}(-\log t)
^{1/2}\,\mathrm{d}
X_{t}^{(\mu)} \biggr) , \qquad  \mu\in I(\mathbb{R}%
^{d}).
\]
Then $\Psi_{-2,2}$ is one-to-one and
%
\begin{equation}
G(\mathbb{R}^{d})=\Psi_{-2,2}(A(\mathbb{R}^{d}))=\Psi_{-2,2}(\Phi
_{\cos
}(I(\mathbb{R}^{d}))). \label{noncom1-2}\vspace*{-2pt}%
\end{equation}
\end{theo}

\begin{pf} Suppose that
$\widetilde{\mu}\in G(\mathbb{R}^{d})$ with triplet $(\widetilde
{\Sigma
},\widetilde{\nu},\widetilde{\gamma})$. Then it follows from Theorems~\ref{GenTGRep} and
\ref{noncom} that
\[
\widetilde{\nu}=\mathcal{A}_{1}(\Upsilon^{0}(\rho))=\Upsilon_{-2.2}
(\mathcal{A}_{1}(\rho))
\]
for some $\rho\in\mathfrak{M}_{L}^{1}(\mathbb{R}^{d})$. Let $\nu
=\mathcal{A}_{1}(\rho)$. Since $\widetilde{\nu}=\Upsilon
_{-2,2}(\nu)$
and since $u=f(t)=(-\log t)^{1/2}$ is the inverse\vadjust{\goodbreak} function of
$t=\int_{u}^{\infty}2v\mathrm{e}^{-v^{2}%
}\,\mathrm{d}v=\mathrm{e}^{-u^{2}}$, we have
\eqref{stochint3} for this $f(t)$
and $T=1$. Choose $\Sigma$ and $\gamma$ satisfying \eqref{stochint2} and
\eqref{stochint4}. Let $\mu\in I(\mathbb{R}^{d})$ having triplet
$(\Sigma
,\nu,\gamma)$. Then $\mu\in A(\mathbb{R}^{d})$ and $\widetilde{\mu
}%
=\Psi_{-2,2}(\mu)$. Notice that $\mathfrak{D}(\Psi_{-2,2})=
I(\mathbb{R}^{d})$, as $f(t)$ is square-integrable on $(0,1)$.
Conversely, we can see that if $\mu\in A(\mathbb{R}^{d})$,
then $\Psi_{-2,2}(\mu)\in G(\mathbb{R}^{d})$. Thus the first
equality in
\eqref{noncom1-2} is proved. The second equality follows from
\eqref{StIntRep2} of Theorem \ref{StIntRep}. The one-to-one property of
$\Psi_{-2,2}$ follows from that of~$\Upsilon_{-2,2}$ in Proposition
\ref{dom-2,2}.
\end{pf}

\begin{rem}
(a) The two representations of $\widetilde{\mu}\in G(\mathbb
{R}^{d})$ in
Theorems \ref{GenTGRep} and \ref{noncom1} are related in the
following way.
Theorem \ref{noncom1} shows that $\widetilde{\mu}\in G(\mathbb
{R}^{d})$ if and
only if $\widetilde{\mu}=\Upsilon_{-2.2}(\Phi_{\cos}(\mu))$ for
some $\mu\in I(\mathbb{R}^{d})$. This $\mu$ has L\'{e}vy measure $\rho^{(1/2)}$
if $\rho$
is the L\'{e}vy measure in the representation of $\widetilde{\mu}$ in Theorem
\ref{GenTGRep}(c). For the proof, use Proposition~\ref{conn2},
Theorems~\ref{StIntRep} and~\ref{noncom}.

(b) We have another representation of the class $G(\mathbb{R}^{d}).$
Let $h(u)=\int_{u}^{\infty
}\mathrm{e}^{-v^{2}}\,\mathrm{d}v$, $u>0$, and denote its inverse
function by
$h^{\ast}(t)$. For $\mu\in I(\mathbb{R}^{d})$, we define
\[
\mathcal{G}(\mu)=\mathcal{L} \biggl( \int_{0}^{\sqrt{\uppi
}/2}h^{\ast
}(t)\,\mathrm{d}X_{t}^{(\mu)} \biggr) .
\]
It is known that $G(\mathbb{R}^{d})=\mathcal{G} ( I(\mathbb{R}%
^{d}) ) $, see Theorem 2.4 (5) in \cite{MS08}. This suggests that
$\mathcal{G}$ is decomposed into
\[
\mathcal{G=}\Psi_{-2.2}\circ\Phi_{\cos}=\Phi_{\cos}\circ\Psi_{-2.2}
\]
with the same domain $I(\mathbb{R}^{d})$, where $\circ$ means
composition of
mappings. The proof is easy to obtain.
\end{rem}


\subsection{$\mathcal{A}_{1}$ is not an Upsilon transformation}\label{sec3.8}


From Theorem \ref{noncom} we obtain the following result.

\begin{theo}
\label{noncom3} The transformation $\mathcal{A}_{1}$ is not an Upsilon
transformation $\Upsilon_{\tau}$ for any dilation measure $\tau$.
\end{theo}

\begin{pf} Suppose that there is a measure $\tau$ on
$(0,\infty)$ such that $\mathcal{A}_{1}(\rho)(B)=\Upsilon_{\tau}
(\rho)(B)$ for $B\in\mathcal{B}(\mathbb{R}^{d})$.
Then, we can show that $\mathcal{A}_{1}(\Upsilon^{0}(\rho))=
\Upsilon^{0}(\mathcal{A}_{1}(\rho))$ for $\rho\in\mathfrak{M}%
_{L}^{1}(\mathbb{R}^{d})$, using the Fubini theorem.
Then, it follows from Theorem \ref{noncom} that
$\Upsilon_{-2,2}(\mathcal{A}_{1}(\rho))=\Upsilon^{0}(\mathcal{A}_{1}
(\rho))$ for $\rho\in\mathfrak{M}_{L}^{1}(\mathbb{R}^{d})$.
Let $\widetilde{\rho}=\mathcal{A}_{1}(\rho).$ If $\int_{\mathbb
{R}^{d}%
} \vert x \vert\rho(\mathrm{d}x)<\infty$, then
\[
\int_{\mathbb{R}^{d}}x\Upsilon^{0}(\widetilde{\rho})(\mathrm{d}x)
=\int%
_{0}^{\infty}\mathrm{e}^{-u}\,\mathrm{d}u\int_{\mathbb{R}^{d}}ux
\widetilde{\rho
}(\mathrm{d}x)=\int_{\mathbb{R}^{d}}x\widetilde{\rho}(\mathrm{d}x),\
\]
while
\[
\int_{\mathbb{R}^{d}}x\Upsilon_{-2,2}(\widetilde{\rho}) (\mathrm
{d}%
x)=\int_{0}^{\infty}2u\mathrm{e}^{-u^{2}}\,\mathrm{d}u \int_{\mathbb
{R}^{d}%
}ux\widetilde{\rho}(\mathrm{d}x)
=2^{-1}\uppi^{1/2} \int%
_{\mathbb{R}^{d}}x\widetilde{\rho}(\mathrm{d}x).
\]
Hence $\Upsilon_{-2,2}(\widetilde{\rho})\neq\Upsilon
^{0}(\widetilde{\rho})$
whenever $\int_{\mathbb{R}^{d}}x\widetilde{\rho}(\mathrm{d}x)\neq
0$ (e.g., choose $\rho=\delta_{e_{1}},e_{1}=(1,0,\ldots,0)$)$.$ This is a
contradiction. Hence the measure $\tau$ does not
exist.\vadjust{\goodbreak}~%
\end{pf}


\section{Examples}\label{sec4}\vspace*{-2pt}


We conclude this paper with examples for Theorems \ref{GenTGRep} and
\ref{noncom}, where the modified Bessel function $K_{0}$ plays a
role. We only
consider the one-dimensional case of L\'{e}vy measures concentrated on
$(0,\infty).$ Multivariate extensions are possible by using polar
decomposition.

By the well-known formula for the modified Bessel functions we have
\[
K_{0}(x)=\frac{1}{2}\int_{0}^{\infty}\mathrm{e}^{-t-x^{2}/(4t)}t^{-1}
\,\mathrm{d}t, \qquad  x>0.\label{exmp1}
\]
An alternative expression is
%
\begin{equation}
K_{0}(x)=\int_{1}^{\infty}(t^{2}-1)^{-1/2}\mathrm{e}^{-xt} \,\mathrm
{d}t, \qquad  x>0; \label{exmp2}%
\end{equation}
see (3.387.3) in \cite{GR07}, page 350. It follows that $K_{0}(x)$ is completely
monotone on $(0,\infty)$ and that $\int_{0}^{\infty}K_{0}(x)\,\mathrm
{d}%
x=\uppi/2$. The Laplace transform of $K_{0}$ in $x>0$ is
%
\begin{equation}
\label{K0}\varphi_{K_{0}}(s):=\int_{0}^{\infty}\mathrm{e}^{-sx}K_{0}
(x)\,\mathrm{d}x=%
\cases{\displaystyle
(1-s^{2})^{-1/2}{\arccos(s)}, &\quad $0<s<1$,\cr\displaystyle
1, &\quad $s=1$,\cr\displaystyle
(1-s^{2})^{-1/2}\log\bigl(s+(s^{2}-1)^{1/2}\bigr), &\quad $s>1$,
}
\end{equation}
see (6.611.9) in \cite{GR07}, page 695.

The following is an example of $\nu$ and $\widetilde{\nu}$ in Theorem
\ref{GenTGRep}(b).\vspace*{-2pt}

\begin{ex}
\label{ex1} Let
\[
\widetilde{\nu}(\mathrm{d}x)=K_{0}(x)1_{(0,\infty)}(x)\,\mathrm{d}%
x\label{exmp3a}
\]
and
%
\begin{equation}
\nu(\mathrm{d}x)=4^{-1}{\uppi}{x^{-1/2}}{\mathrm
{e}^{-x^{1/2}}}1_{(0,\infty)}(x)
\,\mathrm{d}x. \label{exmp4}
\end{equation}
Then $\nu\in\mathfrak{M}_{L}^{B}(\mathbb{R})\cap\mathfrak{M}_{L}%
^{1}(\mathbb{R})$ and $\widetilde{\nu}=\mathcal{A}_{1}(\nu)\in
\mathfrak{M}%
_{L}^{G}(\mathbb{R})$.\vspace*{-2pt}
\end{ex}

The proof is as follows. Since the function ${x^{-1/2}{\mathrm
{e}^{-x^{1/2}}}%
}$ is completely monotone on $(0,\infty)$ and $\int_{0}^{1}x\nu
(\mathrm{d}%
x)<\infty$, we have $\nu\in\mathfrak{M}_{L}^{B}(\mathbb{R})\cap
\mathfrak{M}%
_{L}^{1}(\mathbb{R})$. We can show
%
\begin{equation}\label{new}
\mathcal{A}_{1}(\nu)(B)=\int_{0}^{1}\nu^{(1/2)}(u^{-1}B)2\uppi^{-1}
(1-u^{2})^{-1/2}\,\mathrm{d}u,  \qquad   B\in\mathcal{B}(\mathbb{R}),
\end{equation}
like \eqref{PrUpsIdent1}. Hence
\begin{eqnarray*}
\mathcal{A}_{1}(\nu)(B)
& =&\int_{0}^{1}2^{-1}(1-u^{2})^{-1/2}\,\mathrm{d}u\int_{0}^{\infty}
1_{B}(us^{1/2}){s^{-1/2}}{\mathrm{e}^{-s^{1/2}}}\,\mathrm{d}s\\
& =&\int_{0}^{1}(1-u^{2})^{-1/2}\,\mathrm{d}u\int_{0}^{\infty}1_{B}
(r)\mathrm{e}^{-r/u}\,\mathrm{dr}\\
& =&\int_{0}^{\infty}1_{B}(r)\,\mathrm{d}r\int_{1}^{\infty}(y^{2}-1)^{-1/2}
\mathrm{e}^{-ry}\,\mathrm{d}y\\
& =&\int_{0}^{\infty}1_{B}(r)K_{0}(r)\,\mathrm{d}r=\widetilde{\nu}(B).
\end{eqnarray*}
The fact that $\widetilde{\nu}\in\mathfrak{M}_{L}^{G}(\mathbb{R})$
can also be
shown directly, since $K_{0}({x}^{1/2})$ is again completely monotone in
$x\in(0,\infty)$.

It follows from $\widetilde{\nu}\in\mathfrak{M}_{L}^{G}(\mathbb
{R})$ that
$\widetilde{\nu}$ is the L\'{e}vy measure of a generalized type $G$
distribution $\widetilde{\mu}$ on $\mathbb{R}$. Using (\ref{K0}),
we find that
this $\widetilde{\mu}$ is supported on $[0,\infty)$ if and only if
it has
Laplace transform
\begin{equation}
\int_{\lbrack0,\infty)}\mathrm{e}^{-sx}\widetilde{\mu}(\mathrm{d}x)
=\exp \{ -\gamma_{0}s+\varphi_{K_{0}}(s)-2^{-1}{\uppi} \}
\nonumber
\end{equation}
for some $\gamma_{0}\geq0.$

\begin{rem}
$\mathcal{A}_{1}(\nu)$ in Example \ref{ex1} actually belongs to a smaller
class $\mathfrak{M}_{L}^{B}(\mathbb{R})$. Therefore, in connection to Theorem
\ref{GenTGRep}, it might be interesting to find a necessary and sufficient
condition on $\nu$ in order that $\widetilde{\mu}\in B(\mathbb
{R}^{d})$. The $\nu$
in Example \ref{ex1} also belongs to a smaller class than $\mathfrak
{M}%
_{L}^{B}(\mathbb{R})\cap\mathfrak{M}_{L}^{1}(\mathbb{R})$. It
belongs to the
class of L\'{e}vy measures of distributions in $\mathfrak{R}(\Psi_{-1/2})$
studied in Theorem 4.2 of \cite{Sa06b}.
\end{rem}

  We now give an example of $\rho$ in Theorem \ref{GenTGRep}(c).

\begin{ex}
\label{ex2} Consider the following L\'{e}vy measure in $\mathfrak
{M}_{L}%
^{B}(\mathbb{R})$:
%
\begin{equation}
\rho(\mathrm{d}x)=4^{-1}{{\uppi}^{1/2}}{x^{-1/2}}{\mathrm{e}^{-x/4}}
1_{(0,\infty)}(x)\,\mathrm{d}x. \label{exmp5}%
\end{equation}
Then $\nu$ in (\ref{exmp4}) satisfies $\nu=\Upsilon^{0}(\rho).$
\end{ex}

To prove this, we write $\Upsilon^{0}(\rho)$ as
\[
\Upsilon^{0}(\rho)(\mathrm{d}x)=\int_{0}^{\infty}\rho(u^{-1}
\,\mathrm{d}x)\mathrm{e}^{-u}\,\mathrm{d}u
=4^{-1}{{\uppi}^{1/2}}{x^{-1/2}}\,\mathrm{d}x\int_{0}^{\infty}{u^{-1/2}}
\mathrm{e}^{-u-x/(4u)}\,\mathrm{d}u
\]
on $(0,\infty)$. By formula (3.475.15) in \cite{GR07}, page 369, we have
\[
\int_{0}^{\infty}{u^{-1/2}}\mathrm{e}^{-u-x/(4u)} \,\mathrm{d}u={\uppi}
^{1/2}\mathrm{e}^{-x^{1/2}}.
\]
Hence, $\Upsilon^{0}(\rho)=\nu$ from (\ref{exmp4}).

Since $\mathcal{A}_{1}(\nu)= \mathcal{A}_{1}(\Upsilon^{0}(\rho
))=\Upsilon_{-2,2}(\mathcal{A}_{1}(\rho))$ by Theorem \ref{noncom},
$\mathcal{A}_{1}(\rho)$ is also of interest.

\begin{ex}
Let $\rho$ be as in (\ref{exmp5}). Then
%
\begin{equation}
\mathcal{A}_{1}(\rho)(\mathrm{d}x)={2^{-1}{\uppi}^{-1/2}}\mathrm
{e}^{-x^{2}%
/8}K_{0}(x^{2}/8)1_{(0,\infty)}(x)\,\mathrm{d}x.\label{exmp6}%
\end{equation}
\end{ex}

The proof is as follows. We have, from \eqref{new},
\begin{eqnarray*}
\mathcal{A}_{1}(\rho)(B) & =&{2^{-1}{\uppi}^{-1/2}}\int_{0}^{1}(1-u^{2})^{-1/2}
\,\mathrm{d}u\int%
_{0}^{\infty} 1_{u^{-1}B}(s^{1/2}){s^{-1/2}}{\mathrm{e}^{-s/4}}
\,\mathrm{d}s\\
& =&{{\uppi}^{-1/2}}\int_{0}^{\infty}1_{B}(r)\,\mathrm{d}r\int_{0}^{1}
u^{-1}(1-u^{2})^{-1/2}\mathrm{e}^{-r^{2}/(4u^{2})}\,\mathrm{d}u\\
& =&{2^{-1}{\uppi}^{-1/2}}\int_{0}^{\infty}1_{B}(r)\,\mathrm{d}r \int
_{1}^{\infty
}y^{-1/2}(y-1)^{-1/2}\mathrm{e}^{-r^{2}y/4}\,\mathrm{d}y.
\end{eqnarray*}
Use (3.383.3) in \cite{GR07}, page 347, to obtain
\[
\int_{1}^{\infty}y^{-1/2}(y-1)^{-1/2}\mathrm{e}^{-r^{2}y/4}\,\mathrm
{d}%
y=\mathrm{e}^{-r^{2}/8}K_{0}(r^{2}/8).
\]
Thus we obtain (\ref{exmp6}).

\begin{rem}
The $\rho$ in (\ref{exmp5}) also belongs to $\mathfrak{M}_{L}^{B}
(\mathbb{R})\cap\mathfrak{M}_{L}^{1}(\mathbb{R})$. Therefore,
$\mathcal{A}_{1}
(\rho)$ itself is another example of a measure in
$\mathfrak{M}_L^G(\mathbb{R})$.
\end{rem}

\section*{Acknowledgements}
The authors would like to thank a referee for his/her helpful comments,
which led to an
improvement of the readability of this paper.
Most of this work was done while the second author
visited Keio University in Japan. He gratefully acknowledges the hospitality
and financial support during his stay.

%

\printhistory

\end{document}